# A New Delay Dependent Robust Fixed-Order Decentralized Output Feedback LFC Design for Interconnected Multi-Area Power Systems

Ali Azarbahram, Mahdi Sojoodi[*], Mahmoud-Reza Haghifam

*School of Electrical and Computer Engineering, Tarbiat Modares University, Tehran, Iran.*

**Abstract**:
This paper presents a novel load frequency control (LFC) design using integral-based decentralize fixed-order perturbed dynamic output tracking scheme in a delay dependent nonlinear interconnected multi-area power system via LMI approach. The tracking controller is designed such that the effects of the parameter variations of the plant and dynamic output controller as well as load disturbances are minimized. The inputs and outputs of each controller are local, and these independent controllers are designed such that the robust stability of the overall closed-loop system is guaranteed. Simulation results for both two- and three-area power systems are provided to verify the effectiveness of the proposed design scheme. The simulation results show that the decentralized controlled system behaves well even when there are large parameter perturbations and unpredictable disturbances on the power system.

## 1.     Introduction

Large-scale power systems as geographically distributed systems are composed of several subsystems, i.e., areas including local generator or a group of generators. In addition to its own load, each area is responsible to interchange scheduled power with adjacent areas. The communication between the areas is carried out through tie-lines or high voltage direct current links. A general controller for such system is necessary for keeping the system frequency and the inter-area tie power as near as possible to the scheduled values during the normal operation. Furthermore, in the case of disturbances or sudden changes in load demands, the system should return to its stable situation and scheduled value as soon as possible. In addition to the main goals of using controller, a well-designed power system should also provide high level of power quality. Since all these objectives should be achieved in the presence of uncertainties in modeling, system nonlinearities, and disturbances of the area load, load frequency control (LFC) synthesis is considered as a multi-objective optimization problem [1–3].

LFC problem is discussed from different viewpoints in literature. In [4], a simple robust LFC is introduced based on the Riccati equation approach requiring the bounds of the system parameters. In [5], a PID tuning method for LFC utilizing a two degree of freedom internal model control (IMC) is introduced. In [6], a robust LFC using genetic algorithms and linear matrix inequalities is proposed. A genetic algorithm based fuzzy gain scheduling approach for load frequency control is proposed in [7]. In this approach, a fuzzy system is used to adaptively decide the integral or PI controller gain according to the area control errors (ACE) and their changes.  In [8], a new control system based on the fuzzy sliding mode controller is proposed for controlling the load frequency of nonlinear model of a hydropower plant.

---

[*] Corresponding Author, Address: P.O. Box 14115-194, Tehran, Iran, E-mail address:

*sojoodi@modares.ac.ir*  (M. Sojoodi), web page: *http://www.modares.ac.ir/ece/sojoodi* , Tel./Fax: +98 21 8288-3902.

Recently, decentralized design of interconnected systems has attracted much attention and has been used in the control of multi-area systems, such as interconnected spatially-invariant systems, interconnected computer networks and large-scale power systems. Generally, the works in this area focus on the stability of the system, with local feedback loops in the subsystems as the main goal. In other words, the stability of the system is analyzed independent of the changes in the interconnection topology [9]. However, the control of large-scale interconnected systems is recently shown to benefit from feedback that goes beyond decentralized structures [10]. In [11] a simple and computationally efficient decentralized control design method for the LFC of interconnected power systems is introduced based on the use of a reduced-order observer and a PI controller in each area of the power system. In [12], a decentralized robust LFC for interconnected power systems is considered where each local controller is designed independently under some conditions. In [13], the decentralized load-frequency controller design problem is translated into an equivalent problem of decentralized controller design for a multi-input multi-output control system. The enhancement of inter-area mode damping by multiple flexible AC transmission system devices is demonstrated in [14]. In [15], a decentralized controller is designed for the LFC of interconnected power systems. An adaptive decentralized LFC scheme for multi-area power systems is introduced in [16]. In [17], using an approximation of the interconnection variables, a decentralized sub-optimal LFC scheme is designed as a function of the local area state variables and those resulting from the overlapped states. A decentralized LFC for multi-area power systems is studied in [18] by separating the local transfer matrix from the tie-line power flow network. In [19], a new decentralized robust optimal MISO PID controller based on Characteristic Matrix Eigenvalues and Lyapunov method is proposed. [20] presents a robust decentralized proportional-integral (PI) control design as a solution to the LFC in a multi-area power system where the system robustness margin and transient performance are optimized simultaneously to achieve the optimum PI controller parameters. In [21], an interconnected power system with multi-source power generation is proposed for LFC in deregulated power environment where optimal output feedback controller with pragmatic viewpoint using less number of state variables for feedback is utilized. Robust analysis of decentralized load frequency control (LFC) for multi-area power systems is studied in [22].

Due to the fact that delay usually degrades the performance, and even cause instability of a closed-loop LFC scheme [23], a number of works considering delay in LFC problem are proposed recently. [24] investigates delay-dependent stability of LFC emphasizing on multi-area and deregulated environment. In [25], a delay-dependent robust method is proposed for analysis/synthesis of a PID-type LFC scheme considering time delays. [26] investigates the problem of load frequency control design incorporating the effect of using open communication network instead of a dedicated one for the area control error signals. In [27] an extensive literature review on LFC problem in power system is presented. Some papers report successful application of model predictive control (MPC) technique for LFC. In [28], a feasible Cooperation-Based MPC method is used in distributed LFC. However, the method does not handle the problem of system's parameters mismatch, and only discusses the effect of load change. In [29], an LFC design using the decentralized MPC technique is presented such that the effect of the uncertainty due to governor and turbine parameters variation and load disturbance is reduced compared to the related works. In [30], a robust multivariable MPC is proposed for the solution of LFC in a multi-area power system where The proposed control scheme is designed to consider multivariable nature of LFC, system uncertainty, and generation rate constraint, simultaneously. To achieve robustness against system uncertainty and variation of parameters, a LMI based approach is employed.

Considering nature of LFC problem including plant un-modeled dynamics, quick parameter variations, likelihood of faults and unpredictable disturbances, in addition to computational complexity of the adaptive approaches, the robust control scheme is preferred. The existence of uncertainties in the control implementation has also encouraged the researchers to design non-fragile controllers [31]. Furthermore, there has been a great interest in stabilizing systems by low order or fixed-order controllers. However, their syntheses become non-convex in most cases, and the design conditions are usually given in LMIs but along with certain nonlinear constraints such as rank or inverse constraints [32]. There are some results that lead to BMI form conditions [33].

To the best of our knowledge, no result is available on the robust fixed-order decentralized output feedback LFC design for the delay dependent nonlinear interconnected power systems. Among the models for such systems, the model used in [34] is generally the most applicable. This model assumes that the system is composed of linear subsystems coupled by nonlinear terms with exogenous disturbances. Such a model is useful since the local nonlinear subsystems are usually known with acceptable precision allowing suitable linearization while the interconnections among subsystems are largely unknown and generally nonlinear, but with known bounds. The resulting control scheme offers the maximum possible bounds on the system nonlinearities that still guarantee the stability of the overall interconnected system.

In this paper, decentralized LFC of interconnected power systems is considered. We use a comprehensive dynamical model to design a set of multi-objective robust $H_\infty$ non-fragile fixed-order decentralized output feedback controllers to attenuate disturbances, in the presence of the system nonlinearities, delays and uncertainties as well as the controller parameter uncertainties. The robust controller includes an integral term to achieve the desired load frequency tracking. We derive the sufficient conditions in the form of LMIs, which guarantees the desired $H_\infty$ level of disturbance attenuation while the overall system remains stable.

The rest of this paper is organized as follows: Section 2 provides the problem statement and some preliminary backgrounds. In Section 3, the proposed design method is presented and the theoretical results are provided in the form of a theorem and three corollaries. The case study is discussed in Section 4. Finally, some conclusion remarks are given in Section 5.

## 2. The problem statement and preliminaries

A multi-area power system comprises areas that are interconnected by tie-lines. Figure 1 shows a power system with N control areas. By measuring frequency in each control area against the reference values, the frequency mismatches among interconnections are obtained. The LFC system in each control area of an interconnected power system controls the interchange power with the other control areas as well as its local frequency. Therefore, the tie-line power signal must be considered in the LFC dynamic model.

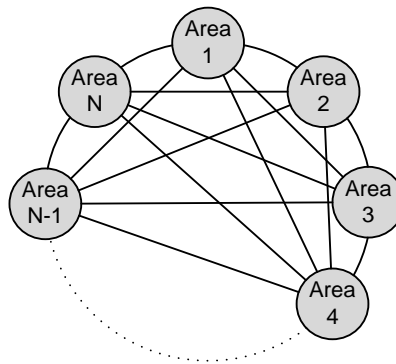

**Figure 1. N-control area power system.**

The nature of the electric power system dynamics is nonlinear in general. Unlike the previous works including [17] and [29] which use linearized models for the LFC problem, we extend the models in [17] and [29] to consider nonlinear terms as in [35], which allows for more accurate controller design for the system. Figure 2 shows the nonlinear block diagram representation of the $i$-th area of a multi-area power system and the basic power systems terminologies are provided in Table 1.

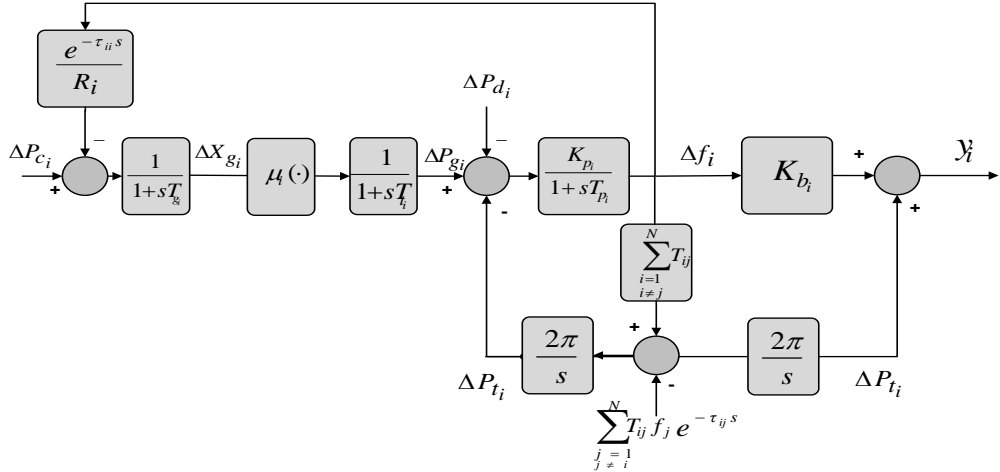

**Figure 2. Nonlinear dynamic model of a control area in an interconnected environment.**

In Figure 2, $T_{ij}$ is the synchronization coefficient between the $i-th$ and the $j-th$ area. $\tau_{ij}$'s are delay terms corresponding to the all areas, where $\tau_{ii}$ is communication delay for the $i-th$ area and $\tau_{ij}$ is communication delay between the $i-th$ and the $j-th$ area. The trajectory of nonlinear function $\mu_i(\Delta X_{gi})$ is depicted in Figure 3.

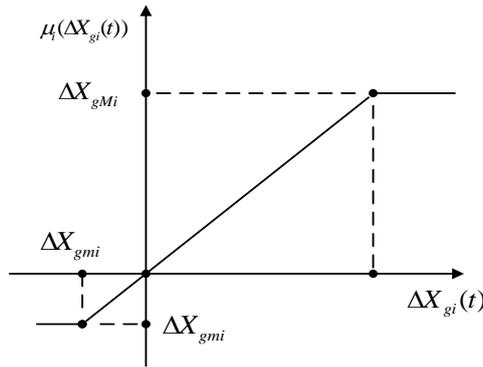

**Figure 3. Nonlinearity of the valve position limits.**

Extending the model in [17] and [29] to include the nonlinear terms, we describe the overall generator-load dynamic equations for the $i-th$ control area with an aggregated generator unit as:

$$\Delta \dot{f}_i(t) = \frac{K_{pi}}{T_{pi}}\Delta P_{gi}(t) - \frac{K_{pi}}{T_{pi}}\Delta P_{di}(t) - \frac{1}{T_{pi}}\Delta f_i(t) - \frac{K_{pi}}{T_{pi}}\Delta P_{ti}(t) \tag{1}$$

$$\Delta \dot{P}_{gi}(t) = \frac{1}{T_{ti}}\eta_i(\Delta X_{gi}(t)) - \frac{1}{T_{ti}}\Delta P_{gi}(t) \tag{2}$$

$$\Delta \dot{X}_{gi}(t) = \frac{1}{T_{gi}}\Delta P_{ci}(t) - \frac{1}{R_i T_{gi}}\Delta f_i(t-\tau_{ii}) - \frac{1}{T_{gi}}\Delta X_{gi}(t) \tag{3}$$

$$\Delta \dot{P}_{ti}(t) = 2\pi\left[\sum_{\substack{j=1\\j\neq i}}^{N}T_{ij}\Delta f_i(t) - \sum_{\substack{j=1\\j\neq i}}^{N}T_{ij}\Delta f_j(t-\tau_{ij})\right]. \tag{4}$$

In a multi-area power system, in addition to regulating the area frequency, a supplementary control maintains the net interchange power with its adjacent areas based on some pre-scheduled values. This is generally accomplished by adding a tie-line flow deviation to the frequency deviation in the supplementary feedback loop. A suitable linear combination of frequency and tie-line power changes for area $i$, is known as the area control error (ACE):

$$y_i(t) = ACE_i = \Delta P_{ti}(t) + K_{b_i}\Delta f_i(t) \tag{5}$$

Assuming $x_{pi}^T = [\Delta X_{gi}^T \quad \Delta P_{gi}^T \quad \Delta f_i^T \quad \Delta P_{ti}^T] \in \Re^4$, $u_i = \Delta P_{ci} \in \Re$, and $w_i = \Delta P_{di} \in \Re$ as the system states, the control input, and the disturbance input, respectively, the nonlinear interconnected system (1)-(5) can be represented as:

$$\dot{x}_{pi}(t) = A_{0ii} x_{pi}(t) + A_{1ii} x_{pi}(t - \tau_{ii}) + \sum_{j=1, j \neq i}^{N} A_{1ij} x_{pj}(t - \tau_{ij}) + B_i u_i(t) + B_{wi} w_i(t) + h_i(t, x_p) \quad (6)$$

$$y_i(t) = C_i x_{pi}(t), \quad i = 1, \cdots, N$$

where $y_i \in \Re$ is the output vector of the $i$-th subsystem and $A_{0ii}$, $A_{1ii}$, $A_{1ij}$, $B_i$, $B_{wi}$, and $C_i$ are known constant matrices with appropriate dimensions, and $h_i(t, x_p)$ is a piecewise continuous vector function in both arguments containing all nonlinearities and uncertainties of system. It is assumed that $h_i(t, x_p)$ satisfies the following quadratic inequalities in its domain of continuity [34–36]:

$$h_i(t, x_p)^T h_i(t, x_p) \leq \overline{\alpha}_i^2 x_p^T \overline{H}_i^T \overline{H}_i x_p \quad i = 1, \ldots, N \quad (7)$$

where $\overline{\alpha}_i > 0$ is a bounding parameter and $\overline{H}_i$ is a constant $v_i \times n$ matrix. Eq. (1)-(5) follows that:

$$A_{0ii} = \begin{bmatrix} -\dfrac{1}{T_{gi}} & 0 & 0 & 0 \\ \dfrac{1}{T_{ti}} & -\dfrac{1}{T_{ti}} & 0 & 0 \\ 0 & \dfrac{K_{pi}}{T_{pi}} & -\dfrac{1}{T_{pi}} & -\dfrac{K_{pi}}{T_{pi}} \\ 0 & 0 & 2\pi \sum_{j=1, j \neq i}^{N} T_{ij} & 0 \end{bmatrix}, A_{1ii} = \begin{bmatrix} 0 & 0 & -\dfrac{1}{R_i T_{gi}} & 0 \\ 0 & 0 & 0 & 0 \\ 0 & 0 & 0 & 0 \\ 0 & 0 & 0 & 0 \end{bmatrix}, A_{1ij} = \begin{bmatrix} 0 & 0 & 0 & 0 \\ 0 & 0 & 0 & 0 \\ 0 & 0 & 0 & 0 \\ 0 & 0 & -2\pi T_{ij} & 0 \end{bmatrix},$$

$$C_i = \begin{bmatrix} 0 & 0 & K_{b_i} & 1 \end{bmatrix}, B_{wi} = \begin{bmatrix} 0 \\ 0 \\ -\dfrac{K_{pi}}{T_{pi}} \\ 0 \end{bmatrix}, h_i(t, x_p) = \begin{bmatrix} 0 \\ \dfrac{1}{T_{ti}} \\ 0 \\ 0 \end{bmatrix} \mu_i(x_{pi}), B_i = \begin{bmatrix} \dfrac{1}{T_{gi}} \\ 0 \\ 0 \\ 0 \end{bmatrix}.$$

Thus, the entire system can be represented as:

$$\dot{x}_p(t) = A_0 x_p(t) + \sum_{i,j=1}^{N} A_{dij} x_p(t - \tau_{ij}) + Bu(t) + B_w w(t) + h(t, x_p) \quad (8)$$

$$y(t) = C x_p(t), \quad i = 1, \cdots, N$$

where, $x_p = \begin{bmatrix} x_{p1}^T & x_{p2}^T & \cdots & x_{pN}^T \end{bmatrix}^T$ and $u(t) = \begin{bmatrix} u_1^T(t) & \cdots & u_N^T(t) \end{bmatrix}^T$, $w = \begin{bmatrix} w_1^T & \cdots & w_N^T \end{bmatrix}^T$, and $y = \begin{bmatrix} y_1^T & \cdots & y_N^T \end{bmatrix}^T$ are the global state vector, input, disturbance, and output vectors, respectively. $A_0 = diag\{A_{011}, A_{022}, \ldots, A_{0NN}\}$, $B = diag\{B_1, B_2, \cdots, B_N\}$, $B_w = diag\{B_{w1}, B_{w2}, \cdots, B_{wN}\}$, $C = diag\{C_1, C_2, \cdots, C_N\}$, and $A_{dij}$ is a block matrix with $N \times N$ blocks where it's $ij-th$ block is $A_{1ij}$ and the other blocks are zero.

$h(t, x_p) = \begin{bmatrix} h_1(t, x_p)^T & \cdots & h_N(t, x_p)^T \end{bmatrix}^T$ is the global nonlinear interconnection vector function. Let $\overline{H}^T = \begin{bmatrix} \overline{H}_1^T & \cdots & \overline{H}_N^T \end{bmatrix}$ where $\overline{H}_i, i = 1, \ldots, N$ are defined in (7), and $\overline{\Gamma}_1 = diag\{\overline{\gamma}_{11} I_{v_1}, \ldots, \overline{\gamma}_{1N} I_{v_N}\}$ with $\overline{\gamma}_{1i} = \overline{\alpha}_i^{-2}$ ($I_{v_i}$ represents the $v_i \times v_i$ identity matrix). Then, it is always possible to find matrices $H$ and $\Gamma_1$ such that:

$$h(t,x_p)^T h(t,x_p) \leq x_p^T \overline{H}^T \overline{\Gamma}_1^{-1} \overline{H} x_p \leq x_p^T H^T \Gamma_1^{-1} H x_p \qquad (9)$$

where $H = diag(H_1, H_2, \cdots, H_N)$, with $v_i \times n_i$ matrices of $H_i$, and:

$$\Gamma_1 = diag\{\gamma_{11} I_{v_1}, ..., \gamma_{1N} I_{v_N}\}, \quad \gamma_{1i} > 0, i = 1, ..., N. \qquad (10)$$

The sufficient condition for (9) gives the following condition for matrices $H$ and $\Gamma_1$ [34]:

$$\lambda_{\max}(\overline{H}^T \overline{H}) \min_i \overline{\gamma}_{1i} \leq \max_i \gamma_{1i} \min_i \lambda_{\min}(H_i^T H_i) \qquad (11)$$

**Assumption 1**: The time-varying delays satisfy:
$$\tau_{ij}(t) > 0,$$
$$|\dot{\tau}_{ij}(t)| \leq d_{ij}, \qquad (12)$$

where $d_{ij}$ is positive known constant.

## 3.  Main results

In this section, a Theorem for fixed-order non-fragile output feedback controller for the delay dependent nonlinear interconnected system given in (8) is presented. The proposed dynamic output feedback controller is the cascade connection of a robust controller and an integrator as shown in Figure 4:

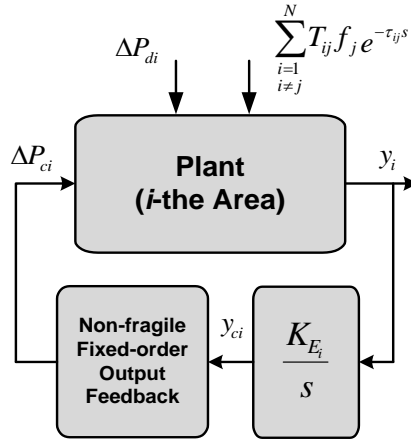

**Figure 4. The proposed controller scheme.**

By augmenting an integral term in each subsystem given in (6), and assuming $x_{pi}^T = [\Delta X_{gi}^T \quad \Delta P_{gi}^T \quad \Delta f_i^T \quad \Delta P_{ti}^T \quad y_{ci}^T] \in \Re^5$, as the augmented states of the $i$-the subsystem, the augmented system parameters which has been described in (6), can now be expressed as:

$$A_{0ii} = \begin{bmatrix} -\frac{1}{T_{gi}} & 0 & 0 & 0 & 0 \\ \frac{1}{T_{ti}} & -\frac{1}{T_{ti}} & 0 & 0 & 0 \\ 0 & \frac{K_{pi}}{T_{pi}} & -\frac{1}{T_{pi}} & -\frac{K_{pi}}{T_{pi}} & 0 \\ 0 & 0 & 2\pi\sum_{\substack{j=1\\j\neq i}}^{N}T_{ij} & 0 & 0 \\ 0 & 0 & K_{Ei}K_{bi} & K_{Ei} & 0 \end{bmatrix}, A_{1ij} = \begin{bmatrix} 0 & 0 & 0 & 0 & 0 \\ 0 & 0 & 0 & 0 & 0 \\ 0 & 0 & 0 & 0 & 0 \\ 0 & 0 & -2\pi T_{ij} & 0 & 0 \\ 0 & 0 & 0 & 0 & 0 \end{bmatrix}, B_i = \begin{bmatrix} \frac{1}{T_{gi}} \\ 0 \\ 0 \\ 0 \\ 0 \end{bmatrix}$$

$$, A_{1ii} = \begin{bmatrix} 0 & 0 & -\frac{1}{R_i T_{gi}} & 0 & 0 \\ 0 & 0 & 0 & 0 & 0 \\ 0 & 0 & 0 & 0 & 0 \\ 0 & 0 & 0 & 0 & 0 \\ 0 & 0 & 0 & 0 & 0 \end{bmatrix}, \begin{aligned} C_i &= \begin{bmatrix} 0 & 0 & 0 & 0 & 1 \end{bmatrix} \\ h_i(t, x_p) &= \begin{bmatrix} 0 & \frac{1}{T_{ti}} & 0 & 0 & 0 \end{bmatrix}^T \mu_i(x_{pi}) \\ B_{wi} &= \begin{bmatrix} 0 & 0 & -\frac{K_{pi}}{T_{pi}} & 0 & 0 \end{bmatrix}^T \end{aligned}$$
(13)

where $K_{E_i}$ is the integral gain for the $i$-th subsystem.

Our objective is to design a non-fragile dynamic output controller for system (8) with the augmented parameters given in (13). Let the non-fragile dynamic output controller be expressed as:

$$C_p : \begin{cases} \dot{x}_c = (A_c + \Delta A_c)x_c + (B_c + \Delta B_c)y_c \\ \Delta P_c = (C_c + \Delta C_c)x_c + (D_c + \Delta D_c)y_c, \ x_c(0) = 0 \end{cases} \tag{14}$$

where $x_c \in \Re^{n_c}$ is the global observer state and $x_c = \begin{bmatrix} x_{c1} & \cdots & x_{cN} \end{bmatrix}^T$, $x_{ci} \in \Re^{n_{ci}}$, $\sum_{i=1}^{N} n_{ci} = n_c$. $A_c = diag\{A_{c1} \cdots A_{cN}\}$, $B_c = diag\{B_{c1} \cdots B_{cN}\}$, $C_c = diag\{C_{c1} \cdots C_{cN}\}$, $D_c = diag\{D_{c1} \cdots D_{cN}\}$, and:

$$\|\Delta A_c\| \leq \delta_{A_c}, \ \|\Delta B_c\| \leq \delta_{B_c}, \ \|\Delta C_c\| \leq \delta_{C_c}, \ \|\Delta D_c\| \leq \delta_{D_c} \tag{15}$$

are norm bounded perturbations. The resulting closed-loop system can be, therefore, represented by:

$$S_{closedp} : \begin{cases} \dot{x}(t) = A_{clp}x(t) + \sum_{i,j=1}^{N} \overline{A}_{dclp-ij}x(t-\tau_{ij}) + E_{clp}w(t) + h^f(t,x) + \begin{bmatrix} B\Delta D_c C x_p + B\Delta C_c x_c \\ \Delta B_c C x_p + \Delta A_c x_c \end{bmatrix} \\ y = C_{cl}x \end{cases} \tag{16}$$

where:

$$x^T = \begin{bmatrix} x_{p1}^T & x_{c1}^T & \cdots & x_{pN}^T & x_{cN}^T \end{bmatrix}, A_{clp} = diag\{A_{clp1} \cdots A_{clpN}\},$$

$$E_{clp} = diag\{E_{clp1} \cdots E_{clpN}\}, h^f(t,x) = \begin{bmatrix} h_1^f(t,x)^T & \cdots & h_N^f(t,x)^T \end{bmatrix}^T$$

$$, h_i^f(t,x) = [h_i(t,x_p)^T \ 0], C_{cl} = diag\{C_{cl1} \cdots C_{clN}\},$$

$$C_{cli} = [C_i \vdots 0], E_{clpi} = \begin{bmatrix} B_{wi} \\ 0 \end{bmatrix}, \overline{A}_{clp-ij} = \begin{bmatrix} A_{1ij} & 0 \\ 0 & 0 \end{bmatrix}$$
(17)

$$A_{clpi} = \begin{bmatrix} A_{0ii} + B_i D_{ci} C_i & B_i C_{ci} \\ B_{ci} C_i & A_{ci} \end{bmatrix}.$$

$\overline{A}_{dclp-ij}$ is also a block matrix with $N \times N$ blocks where it's $ij-th$ block is $\overline{A}_{clp-ij}$ and the other blocks are zero. Also, as a consequence of (9),

$$h^f(t,x)^T h^f(t,x) \leq x^T \overline{H}^{fT} \Gamma_1^{-1} \overline{H}^f x \leq x^T \tilde{H}^T \Gamma_1^{-1} \tilde{H} x_p \tag{18}$$

where $\overline{H}^{fT} = [\overline{H}_1^{fT} \cdots \overline{H}_N^{fT}]$, $\overline{H}_i^f = [\overline{H}_i^1 \ 0 \ \overline{H}_i^2 \ 0 \ ... \ \overline{H}_i^N \ 0]$ in which $v_i \times n_j$ matrices $\overline{H}_i^j (i,j = 1,...,N)$ follow from the decomposition $\overline{H}_i = [\overline{H}_i^1 \cdots \overline{H}_i^N]$ while $\tilde{H} = diag\{\tilde{H}_1,...,\tilde{H}_N\}$ with $\tilde{H}_i = [H_i \ 0]$.

The following two definitions are the vector extensions of those given in [34]:

**Definition 1.** $H_\alpha$ **class**: For any given $\alpha = [\alpha_1^T \cdots \alpha_N^T]^T$ and matrix $H$ with proper dimension, $H_\alpha$ is a class of piecewise-continues functions defined by:

$$H_\alpha = \{h(t,x_p) | h \in \Re^n, h^T h \leq x_p^T H^T \alpha^T \alpha H x_p \text{ in the domains of continuity}\} \tag{19}$$

Therefore, if $h \in H_\alpha$, then $h(t,0) = 0$, which implies that $x_p = 0$ is the equilibrium of system (16).

**Definition 2. Robust stability with vector degree** $\alpha = [\alpha_1^T \cdots \alpha_N^T]^T$: System (16) with bounded nonlinear uncertainty given in (9) is called robustly stable with vector degree $\alpha$ if the equilibrium $x_p = 0$ is globally asymptotically stable for all $h(t,x_p) \in H_\alpha$.

The following theorem gives the sufficient conditions to design a delay dependent integral-based robust non-fragile decentralized fixed-order dynamic output tracker for a nonlinear multi-area power system which ensures $H_\infty$ performance of the closed-loop system (16). Notice that the orders of decentralized controllers can be chosen arbitrarily according to the system conditions and the design constraints and are not restricted by any conditions.

**Theorem 1:** Consider open-loop system (8) and the perturbed control signal given in (14). Given desired weighting scalars $\rho_i$, if positive definite matrix $P = P^T$ in the form of:

$$P = diag(P_1,...,P_N), \ P_i = diag(P_{pi}, P_{ci}), \ i = 1,...,N, \tag{20}$$

and positive definite matrices $R_{ij}$, $\Gamma_1 = diag\{\gamma_{11} I_{v_1},...,\gamma_{1N} I_{v_N}\}$, $\Gamma_2 = diag\{\gamma_{21},...,\gamma_{2N}\}$, and matrices $Y_i$, $Q_i$, $W_i$, $U_i$, for $i = 1,...,N$, $j = 1,...,N$ and positive scalars $s_i, i = 1,..,5$ exist such that the following minimization problem becomes feasible:

minimize $\sum_{i=1}^{N} \gamma_{1i} + \rho_i \cdot \gamma_{2i}$

subject to $\Pi = \begin{bmatrix} \Pi_{11} & \Pi_{12} \\ * & \Pi_{22} \end{bmatrix} \leq 0$

$\Pi_{11} = diag(L_1,...,L_N) + \sum_{i,j=1}^{N} R_{ij}$ (21)

$\Pi_{12} = [P\overline{A}_{dclp-11} \cdots P\overline{A}_{dclp-NN} \ PE_{clp} \ P \ P\begin{bmatrix}B\\0\end{bmatrix} \ P\begin{bmatrix}B\\0\end{bmatrix} \ P\begin{bmatrix}0\\I\end{bmatrix} \ P\begin{bmatrix}0\\I\end{bmatrix} \ s_1 \tilde{H}^T]$

$\Pi_{22} = diag(-(1-d_{11})R_{11} \ ... \ -(1-d_{NN})R_{NN} \ -\Gamma_2 \ -s_1 I \ -s_2 I \ -s_3 I \ -s_4 I \ -s_5 I \ -\Gamma_1)$

where:

$$L_i = \begin{bmatrix} A_{0ii}^T P_{pi} + C_i^T Y_i^T + P_{pi} A_{0ii} + Y_i C_i + C_i^T C_i & C_i^T W_i^T + Q_i \\ Q_i^T + W_i C_i & U_i^T + U_i \end{bmatrix}$$

$$+ \begin{bmatrix} s_3 \|C\| \delta_{D_c}^2 I + s_5 \|C\| \delta_{B_c}^2 I & 0 \\ * & s_2 \delta_{C_c}^2 I + s_4 \delta_{A_c}^2 I \end{bmatrix}, i = 1, \ldots, N. \tag{22}$$

then, the following output feedback controller parameters:

$$A_{ci} = P_{ci}^{-1} U_i, \quad B_{ci} = P_{ci}^{-1} W_i, \quad C_{ci} = B_i^\uparrow P_{pi}^{-1} Q_i, \quad D_{ci} = B_i^\uparrow P_{pi}^{-1} Y_i, \tag{23}$$

with $B_i^\uparrow$ indicating the pseudo inverse of $B_i$, make the overall closed-loop system robustly stable with the robustness degree vector $\begin{bmatrix} 1/\sqrt{\gamma_{11}} & \cdots & 1/\sqrt{\gamma_{1N}} \end{bmatrix}^T$ and simultaneously attenuates disturbances on the output on the $i$-th subsystem to level $\sqrt{\gamma_{2i}}$.

**Proof**: Consider a Lyapunov-Krasovskii functional candidate of the form:

$$V(x_p(t)) = x_p^T P x_p + \sum_{i,j=1}^{N} \int_{t-\tau_{ij}(t)}^{t} x_p^T(\beta) R_{ij} x_p(\beta) d\beta, \quad P > 0, \quad R_{ij} > 0, i, j = 1, \ldots, N., \tag{24}$$

Assuming:

$$z_1 = \Delta C_c x_c, \quad z_2 = \Delta D_c C x_p, \quad z_3 = \Delta A_c x_c, \text{ and } z_4 = \Delta B_c C x_p, \tag{25}$$

one can write:

$$\begin{aligned}
\dot{V}(x_p(t)) &= \dot{x}_p^T P x_p + x_p^T P \dot{x}_p + \sum_{i,j=1}^{N} x_p^T(t) R_{ij} x_p(t) - \sum_{i,j=1}^{N} (1 - \dot{\tau}_{ij}(t))(x_p^T(t - \tau_{ij}(t)) R_{ij} x_p(t - \tau_{ij}(t))) \\
&= (A_{clp} x(t) + \sum_{i,j=1}^{N} \overline{A}_{dclp-ij} x(t - \tau_{ij}) + E_{clp} w(t) + h^f(t, x) + \begin{bmatrix} B z_1 + B z_2 \\ x_3 + z_4 \end{bmatrix})^T P x \\
&+ x^T P (A_{clp} x(t) + \sum_{i,j=1}^{N} \overline{A}_{dclp-ij} x(t - \tau_{ij}) + E_{clp} w(t) + h^f(t, x) + \begin{bmatrix} B z_1 + B z_2 \\ x_3 + z_4 \end{bmatrix}) \\
&+ \sum_{i,j=1}^{N} x_p^T(t) R_{ij} x_p(t) - \sum_{i,j=1}^{N} (1 - \dot{\tau}_{ij}(t))(x_p^T(t - \tau_{ij}(t)) R_{ij} x_p(t - \tau_{ij}(t))).
\end{aligned} \tag{26}$$

Assuming $\Psi^T = \begin{bmatrix} x^T(t) & x^T(t - \tau_{11}(t)) & \cdots & x^T(t - \tau_{NN}(t)) & w^T & h^{fT} & z_1^T & z_2^T & z_3^T & z_4^T \end{bmatrix}$ one can write:

$$\dot{V}(x_p(t)) \leq \Psi^T \overline{\Pi} \Psi \tag{27}$$

where

$$\overline{\Pi} = \begin{bmatrix} \overline{\Pi}_{11} & \overline{\Pi}_{12} \\ * & \overline{\Pi}_{22} \end{bmatrix}$$

$$\overline{\Pi}_{11} = A_{clp}^T P + P A_{clp} + \sum_{i,j=1}^{N} R_{ij} \tag{28}$$

$$\overline{\Pi}_{12} = [P \overline{A}_{dclp-11} \quad \cdots \quad P \overline{A}_{dclp-NN} \quad P E_{clp} \quad P \quad P \begin{bmatrix} B \\ 0 \end{bmatrix} \quad P \begin{bmatrix} B \\ 0 \end{bmatrix} \quad P \begin{bmatrix} 0 \\ I \end{bmatrix} \quad P \begin{bmatrix} 0 \\ I \end{bmatrix}]$$

$$\overline{\Pi}_{22} = diag(-(1 - d_{11}) R_{11} \quad \cdots \quad -(1 - d_{NN}) R_{NN} \quad 0 \quad 0 \quad 0 \quad 0 \quad 0 \quad 0)$$

for $\dot{V}(x_p(t)) < 0$, we should have $-\Psi^T \overline{\Pi} \Psi \geq 0$.

Now, considering an $H_\infty$ disturbance attenuation performance index to attenuate disturbance $w(t)$ [37], and assuming zero initial conditions, the following must be satisfied:

$$\dot{V} + y^T(t) y(t) - w^T(t) \Gamma_2 w(t) < 0. \tag{29}$$

Substituting (16) in (29) results in:
$$\dot{V}(t) + x^T(t)C_{clp}^T C_{clp} x(t) - w^T(t)\Gamma_2 w(t) < 0. \quad (30)$$

Using S-Procedure [32], we conclude from (28) and (30) that:
$$\hat{\Pi} = \begin{bmatrix} \hat{\Pi}_{11} & \hat{\Pi}_{12} \\ * & \hat{\Pi}_{22} \end{bmatrix} < 0$$

$$\hat{\Pi}_{11} = A_{clp}^T P + PA_{clp} + \sum_{i,j=1}^{N} R_{ij} + C_{clp}^T C_{clp} \quad (31)$$

$$\hat{\Pi}_{12} = \overline{\Pi}_{12}$$

$$\hat{\Pi}_{22} = diag\left(-(1-d_{11})R_{11} \quad \cdots \quad -(1-d_{NN})R_{NN} \quad -\Gamma_2 \quad 0 \quad 0 \quad 0 \quad 0 \quad 0 \right).$$

From (18) and (25), one can conclude:
$$z_1^T z_1 \leq \delta_{C_c}^2 x_c^T x_c, \quad z_2^T z_2 \leq \|C\|\delta_{D_c}^2 x_p^T x_p, \quad z_3^T z_3 \leq \delta_{A_c}^2 x_c^T x_c, \quad z_4^T z_4 \leq \|C\|\delta_{B_c}^2 x_p^T x_p, \quad h^{fT} h^f \leq x^T \tilde{H}^T \Gamma_1^{-1} \tilde{H} x. \quad (32)$$

Assuming (20) and using *S-Procedure* and *Schur-Complement* [32], one can conclude from (31) that:

$$\tilde{\Pi} = \begin{bmatrix} \tilde{\Pi}_{11} & \tilde{\Pi}_{12} \\ * & \tilde{\Pi}_{22} \end{bmatrix} < 0$$

$$\tilde{\Pi}_{11} = diag(L_1, \ldots, L_N) + \sum_{i,j=1}^{N} R_{ij}$$

$$\tilde{\Pi}_{12} = [P\overline{A}_{dclp-11} \quad \cdots \quad P\overline{A}_{dclp-NN} \quad PE_{clp} \quad P \quad P\begin{bmatrix}B\\0\end{bmatrix} \quad P\begin{bmatrix}B\\0\end{bmatrix} \quad P\begin{bmatrix}0\\I\end{bmatrix} \quad P\begin{bmatrix}0\\I\end{bmatrix} \quad s_1 \tilde{H}^T] \quad (33)$$

$$\tilde{\Pi}_{22} = diag(-(1-d_{11})R_{11} \ldots -(1-d_{NN})R_{NN} \quad -\Gamma_2 \quad -s_1 I$$
$$-s_2 I \quad -s_3 I \quad -s_4 I \quad -s_5 I \quad -\Gamma_1)$$

where:
$$L_i = \begin{bmatrix} A_{0ii}^T P_{pi} + C_i^T D_{ci}^T B_i^T P_{pi} + P_{pi} A_{0ii} + P_{pi} B_i D_{ci} C_i + C_i^T C_i & C_i^T B_{ci}^T P_{ci} + P_{pi} B_i C_{ci} \\ C_{ci}^T B_i^T P_{Pi} + P_{ci} B_{ci} C_i & A_{ci}^T P_{ci} + P_{ci} A_{ci} \end{bmatrix}$$
$$+ \begin{bmatrix} s_3 \|C\|\delta_{D_c}^2 I + s_5 \|C\|\delta_{B_c}^2 I & 0 \\ * & s_2 \delta_{C_c}^2 I + s_4 \delta_{A_c}^2 I \end{bmatrix}, i = 1, \ldots, N. \quad (34)$$

Now, changing the variables as:
$$Y_i = P_{pi} B_i D_{ci}, \quad Q_i = P_{pi} B_i C_{ci}, \quad W_i = P_{ci} B_{ci}, \quad U_i = P_{ci} A_{ci}, \quad (35)$$

(21) is obtained. This completes the proof.□

The following corollary gives the simplified result for the special case that the controller parameters can be implemented accurately.

**Corollary 1:** Consider open-loop system (8) and the controller (14) without uncertainty, i.e., $\Delta A_c = 0$, $\Delta B_c = 0$, $\Delta C_c = 0$, $\Delta D_c = 0$. Given desired weighting scalars $\rho_i$, if positive definite matrix $P = P^T$ of the form (20) and positive definite matrices $R_{ij}$, $\Gamma_1 = diag\{\gamma_{11} I_{v_1}, \ldots, \gamma_{1N} I_{v_N}\}$, $\Gamma_2 = diag\{\gamma_{21}, \ldots, \gamma_{2N}\}$ and matrices $Y_i$, $Q_i$, $W_i$, $U_i$ for $i = 1, \ldots, N$, $j = 1, \ldots, N$ and positive scalar $s_1$ exist, such that the following minimization problem is feasible:

$$\text{minimize} \quad \sum_{i=1}^{N} \gamma_{1i} + \rho_i \cdot \gamma_{2i}$$

$$\text{subject to} \quad \Pi = \begin{bmatrix} \Pi_{11} & \Pi_{12} \\ * & \Pi_{22} \end{bmatrix} \leq 0$$

$$\Pi_{11} = diag(L_1,\ldots,L_N) + \sum_{i,j=1}^{N} R_{ij} \tag{36}$$

$$\Pi_{12} = [P\overline{A}_{dclp-11} \quad \cdots \quad P\overline{A}_{dclp-NN} \quad PE_{clp} \quad P \quad s_1\tilde{H}^T]$$

$$\Pi_{22} = diag(-(1-d_{11})R_{11} \quad \cdots \quad -(1-d_{NN})R_{NN} \quad -\Gamma_2 \quad -s_1 I \quad -\Gamma_1)$$

where:

$$L_i = \begin{bmatrix} A_{0ii}^T P_{pi} + C_i^T Y_i^T + P_{pi} A_{0ii} + Y_i C_i + C_i^T C_i & C_i^T W_i^T + Q_i \\ Q_i^T + W_i C_i & U_i^T + U_i \end{bmatrix}, i = 1,\ldots,N. \tag{37}$$

Then, the output feedback controller parameters given in (23) make the overall closed-loop system robustly stable with the robustness degree vector $[1/\sqrt{\gamma_{11}} \quad \cdots \quad 1/\sqrt{\gamma_{1N}}]^T$ and simultaneously attenuate disturbances on the output on the $i-th$ subsystem to level $\sqrt{\gamma_{2i}}$.

**Corollary 2:** Considering nonlinear dynamic model shown in Figure 2, and setting all delay terms $\tau_{ij}'s$ equal to zero, constituting open-loop system and using controller (14), then the closed-loop system is achieved, and the LMI minimization problem can be solved similar to Theorem 1.

**Corollary 3:** Considering Corollary 2 and using controller (14), without uncertainty, i.e., $\Delta A_c = 0$, $\Delta B_c = 0$, $\Delta C_c = 0$, $\Delta D_c = 0$, one can solve an LMI problem in the same way as corollary 1 and the parameters of controller can consequently be achieved.

**Remark 1**: Considering zero as the order of dynamic output controller, the controllers in Corollary 1 and Theorem 1 reduce to static output feedback. In addition, if one assumes $C = I$, then the dynamic output controller design problem changes to state feedback controller design problem.

## 4. Simulation results

The results of Theorem 1 and Corollary 1-Corollary 3 are applied to the two- and three-area power systems given in [17] and [29] respectively, whose models are described by (1)-(5).

The simulation results are presented and discussed in this section in order to validate the effectiveness of the proposed scheme. The Matlab/Simulink software package are used for this purpose.

### 4.1. Two-control area

In this section, two-control area power system is considered as a test system to illustrate effectiveness of the proposed control strategy.

#### 4.1.1. Nominal performance without delay

The nominal parameter values of a practical two-area power system [29], and the parameters of integral based fixed-order output feedback controller derived from the design in Corollary 3, where there is no delay term in power system model and the controller is not perturbed, are given in Table 2. Here, the order of local output feedback controllers are selected as $n_{c_1} = n_{c_2} = 2$; however as stated before, this value can be chosen arbitrarily according to the system conditions and the design constraints. Let weighting parameters be chosen as: $\rho_i = 1$, $(i = 1,2)$.

The system performance with the proposed output feedback controllers at nominal parameters is measured and compared with the system performance with a conventional integrator and MPC proposed in [29] by assuming $\Delta P_{d_1}$, $\Delta P_{d_2}$ being equal to 0.02 $p.u.$ at $t = 3$ sec. Figure 5 (e) and (f) show the frequency deviations of the two areas $(\Delta f_1, \Delta f_2)$ in this case, using the derived fixed order output feedback controller defined in Corollary 3, the proposed MPC in [29], and the conventional integral controller following a step load change in both area, Respectively.

**4.1.2. Robustness validation for the model without delay terms (corollary 2)**

To verify the robustness and the non-fragility of the proposed approach without model delay terms, we increase the controller gains by 10%, and also simultaneously perturb the system parameters as given in Table 3. Figure 5 ((a) to (d)) depict the response of the proposed approach investigated in Corollary 2, MPC and integral controller with the same load change as described in the nominal case. The results show that the proposed controller outperforms the MPC controller in [29] while the conventional integral control method leads to instability under these perturbations.

It is clear that the closed-loop system response with the proposed controller is faster and more stable compared to the system with traditional integrator and MPC controller. It is notable that although the MPC approach in [29] is adaptive, our design is robust. This further highlights the performance of the proposed design approach.

Considering the nonlinearities as Figure 3, the robustness degree of each subsystem must be greater than one. The obtained results for $\Gamma_1$ in Table 2 and Table 3 show that in addition to the nonlinear terms, the system uncertainties are also covered appropriately.

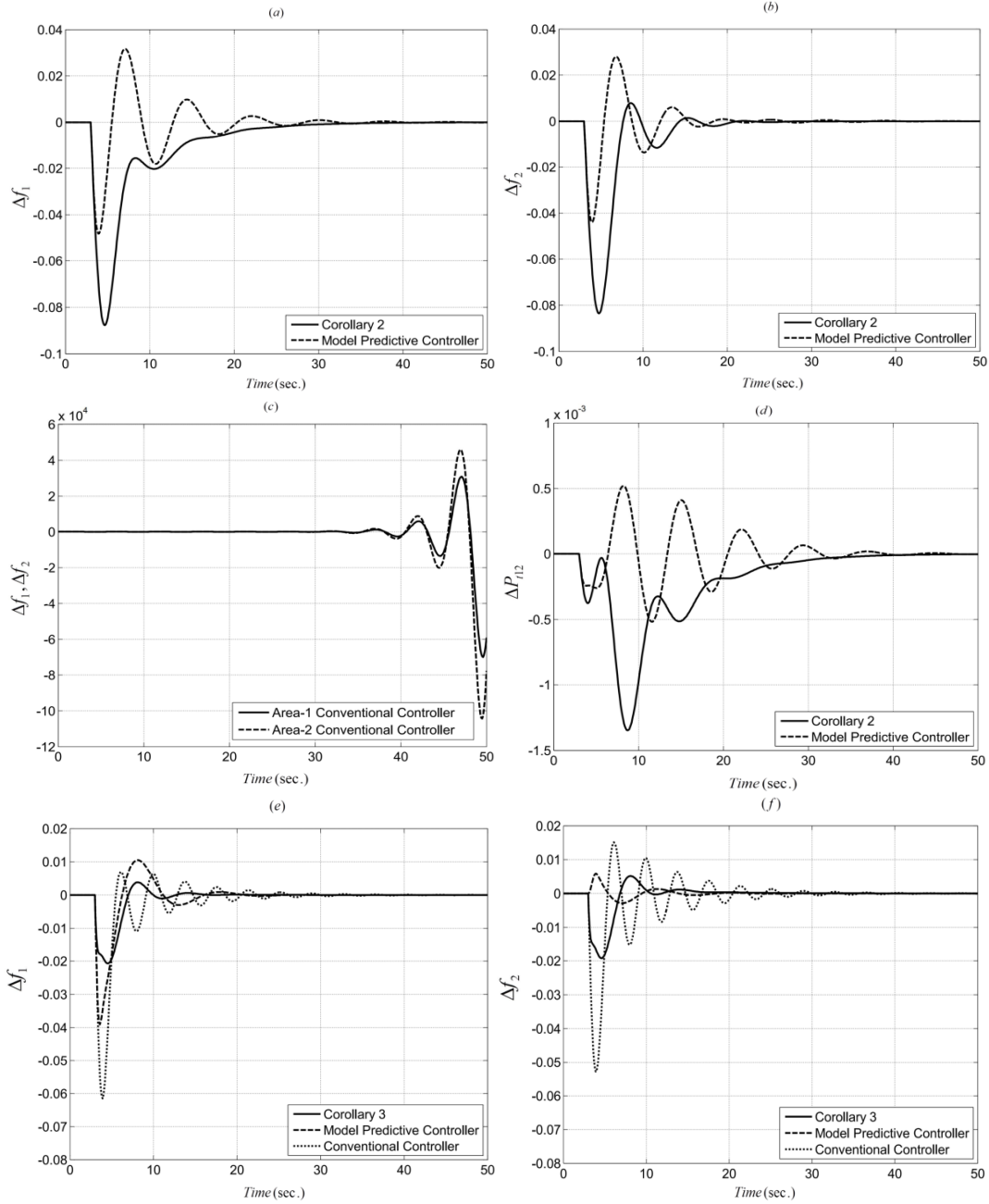

**Figure 5.** (a) Frequency deviation in area-1 (b) Frequency deviation in area-2 for: the proposed controller in Corollary 2 (solid), MPC (dashed) (c) Frequency deviation in area-1 (solid) and area-2 (dashed) with conventional controller (d) Tie-line power change between area-1 and area-2 for: the proposed controller in Corollary 2 (solid), MPC (dashed) (e) Frequency deviation in area-1 (f) Frequency deviation in area-2 for: the proposed controller in Corollary 3 (solid), MPC (dashed).

### 4.1.3. Nominal performance for the model with considering delay

Considering nominal parameters of a practical two-area power system [29] in the same way as Table 2 and assuming all delay terms being equal to $2+0.3\sin t$ seconds and setting $K_{E_i}$ for each subsystem to -1.8, the parameters of delay dependent integral based fixed-order output feedback controller derived from the design in Corollary 1 are given in Table 4. Figure 6 shows the simulation results in this case, where (a) and (b) depict frequency deviation of area-1 and area-2 according to the proposed approach in Corollary 1 (solid) and MPC approach in [29] (dashed), (c) demonstrate the same but for Conventional controller, and (d) shows Tie-line power change between area-1 and area-2.

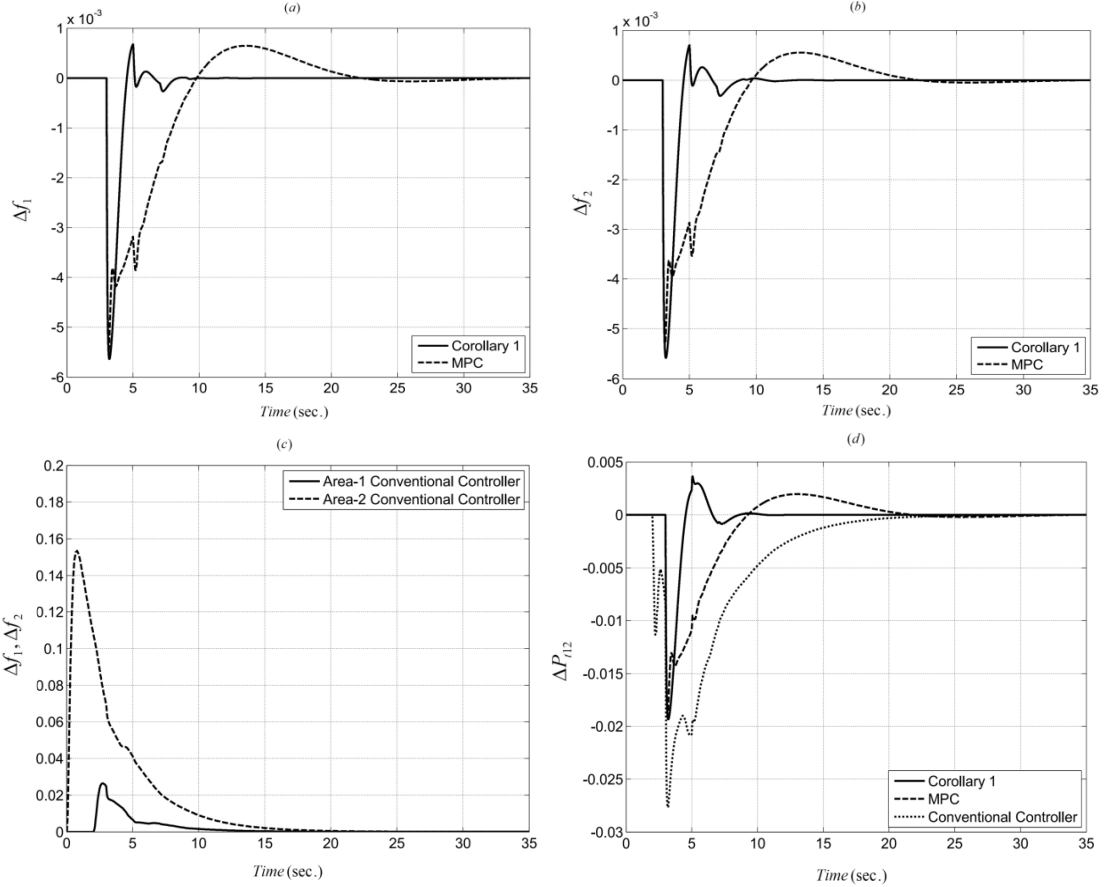

**Figure 6.** (a) Frequency deviation in area-1 (b) Frequency deviation in area-2 for: the proposed controller in Corollary 1 (solid), MPC (dashed) (c) Frequency deviation in area-1 (solid) and area-2 (dashed) with conventional controller (d) Tie-line power change between area-1 and area-2 for: the proposed controller in Corollary 1 (solid), MPC (dashed).

It can be seen from Figure 6 that the closed-loop system response with the proposed controller defined in Corollary 1 is faster and more stable compared to the system with traditional integrator and MPC controller and also, the overshoot percent of traditional controller is partly large compared to the proposed controller defined in Corollary 1 and MPC controller.

**4.1.4. Robustness validation for the model with delay terms**

Again in order to verify the robustness non-fragility of the proposed approach for the model of a delay dependent power system we increase the controller gains by 10%, and simultaneously perturb the system parameters in the same way as given in Table 3. Assuming all delay terms being equal to $2+0.3\sin t$ seconds and setting $K_{E_i}$ for each subsystem to -2.3, the parameters of delay dependent integral based fixed-order output feedback controller derived from the design in Theorem 1 are given in Table 5. Figure 7 depicts the response of the proposed approach in Theorem 1 where (a) and (b) depict frequency deviation of area-1 and area-2, (c) shows Tie-line power change between area-1 and area-2 according to the proposed approach in Theorem 1 (solid) and MPC approach in [29] (dashed), and (d) illustrates frequency deviation of area-1 and area-2 with Conventional controller.

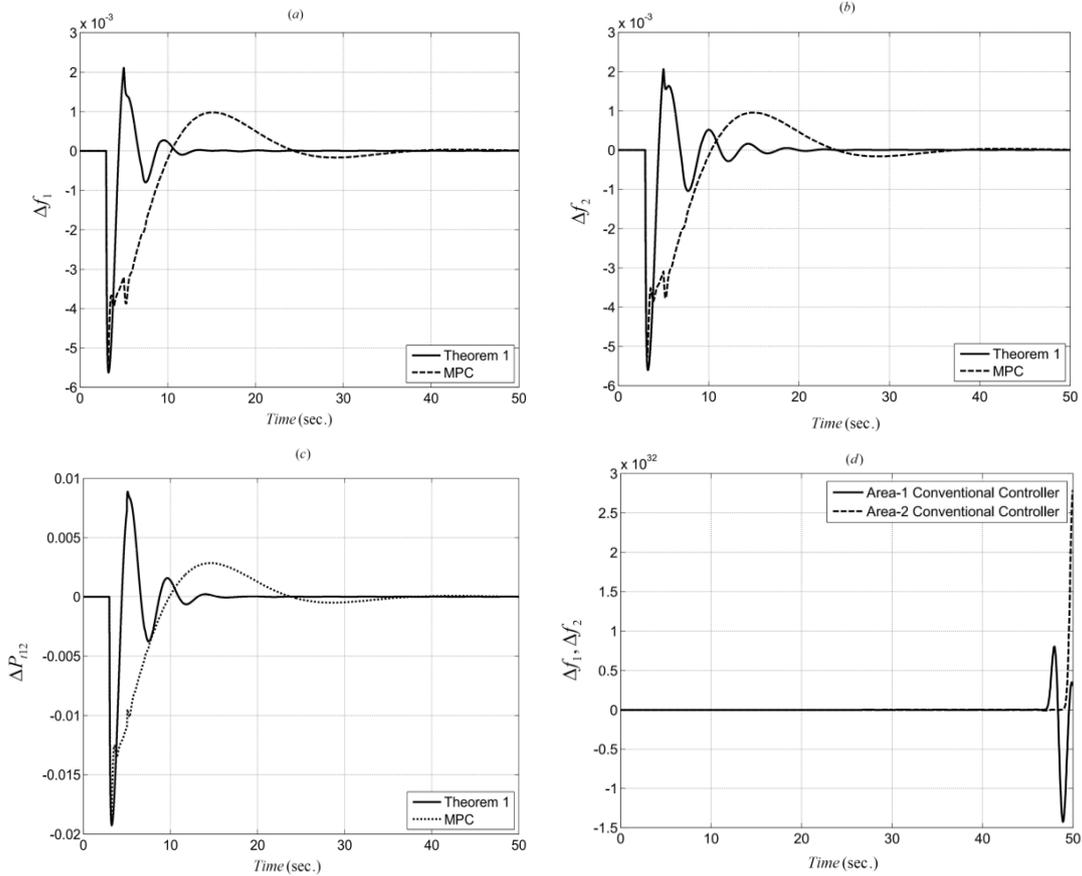

**Figure 7. (a) Frequency deviation in area-1 (b) Frequency deviation in area-2 (c) Tie-line power change between area-1 and area-2 for: the proposed controller in Theorem 1 (solid), MPC (dashed) (d) Frequency deviation in area-1 (solid) and area-2 (dashed) with conventional controller.**

The results show that the proposed controller in Theorem 1 outperforms the MPC controller in [29] while the conventional integral control method leads to instability under these perturbations.

## 4.2. Three-control area

In this sub-section, another multi-area power system with three interconnected control areas is considered as given in [17]. The nominal parameter values of the three-area power system [17] and the derived controller parameters according to Corollary 1 are given in Table 6. Assume that all delay terms are equal to $2+0.3\sin t$ seconds and weighting parameters are chosen as: $\rho_i = 1$, $(i=1,2,3)$.

Figure 8 shows the same results as reported in section 4.1.3 for the case of three area LFC problem. The results show improvement comparing to the proposed design approach with MPC and conventional controller as discussed in section 4.1.3.

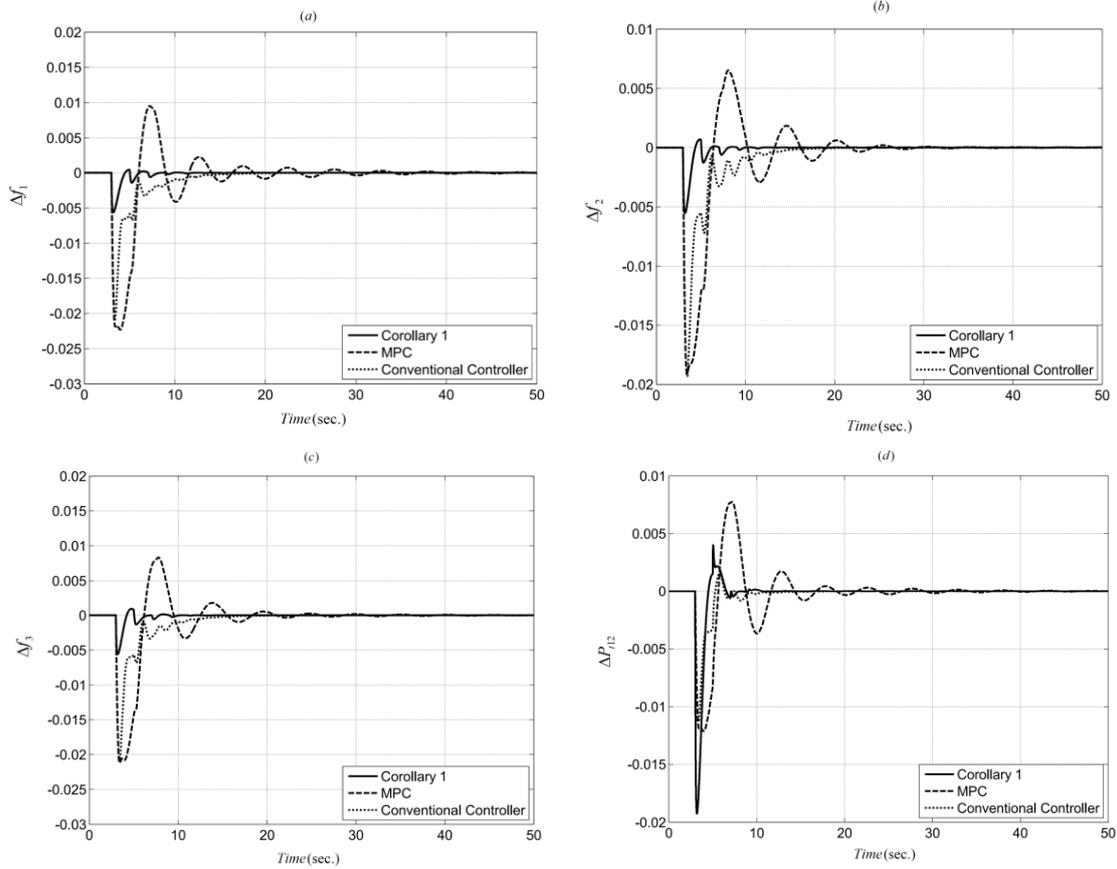

**Figure 8. (a) Frequency deviation in area-1 (b) Frequency deviation in area-2 (c) Frequency deviation in area-3 (d)Tie-line power change between area-1 and area-2 for: the proposed controller in Corollary 1 (solid), MPC (dashed) and conventional (dotted) controller.**

## 5. Conclusion

This paper proposed an integral based decentralized control scheme for the LFC of a delay dependent multi-area power system, treated as an interconnected dynamical system. To achieve this goal, a new LMI-based procedure has been formulated for the design of decentralized fixed-order dynamic output controllers for the systems composed of linear subsystems coupled by uncertain nonlinear interconnections satisfying quadratic constraints. The local Lyapunov-Krasovskii matrices together with the corresponding robustness degrees and prescribed disturbance attenuation levels are derived by solving a convex optimization problem. Then the parameters of the fixed-order controllers for the subsystems are extracted which robustly stabilizes the overall system. The order of controllers can be chosen arbitrarily according to the system conditions and limitations. To verify the proposed control scheme, simulations are carried out on a two and three-area power system. The results of simulations showed that the derived controllers are robust in the presence of controllers' gain perturbations, and the power system disturbances and parameter variations.

**Table 1. Basic power systems terminology of the $i$-th area**

| Parameter | Description |
|---|---|
| $T_{p_i}$ | Plant model time constant |
| $T_{t_i}$ | Turbine time constant |
| $T_{g_i}$ | Governor time constant |
| $K_{p_i}$ | Plant gain |
| $R_i$ | Speed regulation due to governor action |
| $K_{b_i}$ | Frequency biasing factor |
| $\Delta f_i$ | Incremental changes in the frequency |
| $\Delta P_{gi}$ | Incremental changes in the generator output power |
| $\Delta X_{gi}$ | Incremental changes in the governor valve position |
| $\Delta P_{ti}$ | Tie-line power deviation |
| $\Delta P_{ci}$ | Control signal |
| $\Delta P_{di}$ | Load change |
| $y_i$ | Area Control Error |

# Appendix

**Table 2. Nominal parameters of the two-area power system without delay and designed controller parameters**

| Parameter | Area 1 ($i=1$) | Area 2 ($i=2$) |
|---|---|---|
| $T_{p_i}$ | 11.1133 | 12.6 |
| $T_{t_i}$ | 0.4 | 0.44 |
| $T_{g_i}$ | 0.08 | 0.06 |
| $K_{p_i}$ | 66.667 | 62.5 |
| $R_i$ | 3 | 2.73 |
| $K_{E_i}$ | -0.03 | -0.02 |
| $K_{b_i}$ | 1 | 1 |
| $\Delta X_{gmi}$ | -0.03 | -0.03 |
| $\Delta X_{gMi}$ | 0.12 | 0.12 |
| $2\pi \sum_{i,j=1}^{N} T_{ij}$ | 1.2566 | 1.2566 |
| $A_c (n_{c_1} = n_{c_2} = 2)$ | $\begin{bmatrix} -1.2132 & -0.0497 \\ 0.0497 & -1.2132 \end{bmatrix}$ | $\begin{bmatrix} -1.2264 & -0.0632 \\ 0.0632 & -1.2264 \end{bmatrix}$ |
| $B_c$ | $\begin{bmatrix} -0.3645 \\ 0.3645 \end{bmatrix}$ | $\begin{bmatrix} -0.3741 \\ 0.3741 \end{bmatrix}$ |
| $C_c$ | $\begin{bmatrix} 0.9174 & 0.9174 \end{bmatrix}$ | $\begin{bmatrix} 0.9617 & 0.9617 \end{bmatrix}$ |
| $D_c$ | 4.25 | 4.2276 |

| | $\gamma_{1i}$ | 0.38 | 0.38 |
|---|---|---|---|
| | $\gamma_{2i}$ | 0.2773 | 0.2773 |
| | $h_i(t, x_p)$ | $\begin{bmatrix} 0 \\ 1 \\ \frac{1}{T_{t1}} \\ 0 \\ 0 \end{bmatrix} \mu_1(x_{p1})$ | $\begin{bmatrix} 0 \\ 1 \\ \frac{1}{T_{t2}} \\ 0 \\ 0 \end{bmatrix} \mu_2(x_{p2})$ |

**Table 3. The perturbed parameters of the closed-loop two-area power system without delay terms**

| Parameter | Area 1 ($i=1$) | changes | Area 2 ($i=2$) | Changes |
|---|---|---|---|---|
| $T_{p_i}$ | 22.2266 | 100% | 25.2 | 100% |
| $T_{t_i}$ | 0.8 | 100% | 0.88 | 100% |
| $T_{g_i}$ | 0.105 | 31.25% | 0.105 | 75% |
| $K_{p_i}$ | 100 | 50% | 100 | 60% |
| $R_i$ | 3.5 | 16% | 3.5 | 26% |
| $K_{E_i}$ | -.4 | 33.33% | -.3 | 50% |
| $K_{b_i}$ | 1.5 | 50% | 1.5 | 50% |
| $\Delta X_{gmi}$ | -0.03 | 10% | -0.03 | 10% |
| $\Delta X_{gMi}$ | 0.12 | 10% | 0.12 | 10% |
| $2\pi \sum_{i,j=1}^{N} T_{ij}$ | 2.5132 | 100% | 2.5132 | 100% |
| $\delta_{A_c}, \delta_{B_c}, \delta_{C_c}, \delta_{D_c}$ | 0.1 | | 0.1 | |
| $A_c (n_{c_1} = n_{c_2} = 2)$ | $\begin{bmatrix} -1.2465 & -0.0837 \\ -0.0837 & -1.2465 \end{bmatrix}$ | | $\begin{bmatrix} -1.2439 & -0.079 \\ -0.079 & -1.2439 \end{bmatrix}$ | |
| $B_c$ | $\begin{bmatrix} -0.3775 \\ -0.3775 \end{bmatrix}$ | | $\begin{bmatrix} -0.38 \\ -0.38 \end{bmatrix}$ | |
| $C_c$ | $[0.933 \; 0.933]$ | | $[0.954 \; 0.954]$ | |
| $D_c$ | 4.1743 | | 4.2027 | |
| $\gamma_{1i}$ | 0.1197 | | 0.1197 | |
| $\gamma_{2i}$ | 0.24 | | 0.24 | |

**Table 4. Controller parameters of the two-area power system with fixed delay terms**

| Parameter | Area 1 ($i=1$) | Area 2 ($i=2$) |
|---|---|---|
| $A_c (n_{c_1} = n_{c_2} = 2)$ | $\begin{bmatrix} -1.5923 & -1.5897 \\ -1.5897 & -1.5923 \end{bmatrix} * 10^4$ | $\begin{bmatrix} -1.1217 & -1.1176 \\ -1.1176 & -1.1217 \end{bmatrix} * 10^4$ |
| $B_c$ | $\begin{bmatrix} -1.5891 \\ -1.5891 \end{bmatrix} * 10^4$ | $\begin{bmatrix} -1.1183 \\ -1.1183 \end{bmatrix} * 10^4$ |
| $C_c$ | $[0.0526 \; 0.0526] * 10^4$ | $[0.0283 \; 0.0283] * 10^4$ |
| $D_c$ | $0.0525 * 10^4$ | $0.0282 * 10^4$ |
| $\gamma_{1i}$ | 0.1245 | 0.1245 |
| $\gamma_{2i}$ | 0.3224 | 0.3224 |

**Table 5. The parameters of the proposed method with fixed delay terms**

| Parameter | Area 1 ($i=1$) | Area 2 ($i=2$) |
|---|---|---|
| $A_c (n_{c_1} = n_{c_2} = 2)$ | $\begin{bmatrix} -1.59232 & -1.589708 \\ -1.589708 & -1.59232 \end{bmatrix} * 10^4$ | $\begin{bmatrix} -1.12171 & -1.117605 \\ -1.117605 & -1.12171 \end{bmatrix} * 10^4$ |
| $B_c$ | $\begin{bmatrix} -1.589112 \\ -1.589112 \end{bmatrix} * 10^4$ | $\begin{bmatrix} -1.118311 \\ -1.118311 \end{bmatrix} * 10^4$ |
| $C_c$ | $[0.052616 \ 0.052616] * 10^4$ | $[0.028311 \ 0.028311] * 10^4$ |
| $D_c$ | $0.05251 * 10^4$ | $0.028209 * 10^4$ |
| $\gamma_{1i}$ | 0.1427 | 0.1427 |
| $\gamma_{2i}$ | 0.28 | 0.28 |

**Table 6. Nominal parameters of the three-area power system.**

| Parameter | Area 1 ($i=1$) | Area 2 ($i=2$) | Area 3 ($i=3$) |
|---|---|---|---|
| $T_{p_i}$ | 20 | 25 | 20 |
| $T_{t_i}$ | 0.3 | 0.33 | 0.35 |
| $T_{g_i}$ | 0.08 | 0.072 | 0.07 |
| $K_{p_i}$ | 120 | 112.5 | 115 |
| $R_i$ | 2.4 | 2.7 | 2.5 |
| $K_{E_i}$ | -1.85 | -1.85 | -1.85 |
| $K_{b_i}$ | 1 | 1 | 1 |
| $\Delta X_{gmi}$ | -0.03 | -0.03 | -0.03 |
| $\Delta X_{gMi}$ | 0.12 | 0.12 | 0.12 |
| $2\pi \sum_{i,j=1}^{N} T_{ij}$ | 0.545 | 0.545 | 0.545 |
| $A_c (n_{c_1} = n_{c_2} = n_{c_3} = 2)$ | $\begin{bmatrix} -1.5923 & -1.5897 \\ -1.5897 & -1.5923 \end{bmatrix} * 10^4$ | $\begin{bmatrix} -1.1217 & -1.1176 \\ -1.1176 & -1.1217 \end{bmatrix} * 10^4$ | $\begin{bmatrix} -1.0817 & -1.0676 \\ -1.0676 & -1.0817 \end{bmatrix} * 10^4$ |
| $B_c$ | $\begin{bmatrix} -1.5891 \\ -1.5891 \end{bmatrix} * 10^4$ | $\begin{bmatrix} -1.1183 \\ -1.1183 \end{bmatrix} * 10^4$ | $\begin{bmatrix} -1.0583 \\ -1.0583 \end{bmatrix} * 10^4$ |
| $C_c$ | $[0.0526 \ 0.0526] * 10^4$ | $[0.0283 \ 0.0283] * 10^4$ | $[0.0383 \ 0.0383] * 10^4$ |
| $D_c$ | 525 | 282 | 582 |
| $\gamma_{1i}$ | 0.1621 | 0.1621 | 0.1621 |
| $\gamma_{2i}$ | 0.23 | 0.23 | 0.23 |